\documentclass[11pt,a4paper,leqno,draft]{article}
\usepackage{amstex}

\title{An asymptotic existence theorem for plane curves with 
prescribed singularities}
\author{Thierry Mignon\\
{\small \textit{ENS Lyon, UMR CNRS 5669}}\\
{\small \textit{46 All\'ee d'Italie, 69364 Lyon Cedex, France}}\\
{\small \textit{e-mail~: tmignon@@umpa.ens-lyon.fr}}}

\date{}

\newtheorem{theorem}{Theorem}

\newtheorem{proposition}{Proposition}[section]
\newtheorem{lemma}[proposition]{Lemma}
\newtheorem{corollary}[proposition]{Corollary}
\newtheorem{definition}[proposition]{Definition}
\newtheorem{definitions}[proposition]{Definitions}
\newtheorem{remark}[proposition]{Remark}

\newcommand{\al}{\alpha}
\newcommand{\eps}{\varepsilon}
\newcommand{\fl}{\longrightarrow}
\newcommand{\isom}{\stackrel{\sim}{\longrightarrow}}
\newcommand{\equi}{\Longleftrightarrow}
\newcommand{\morph}[1]{\stackrel{#1}{\longrightarrow}}

\newcommand{\p}{\Bbb P}
\newcommand{\Z}{ \Bbb Z}

\newcommand{\sq}{\ensuremath{\square}}
\newcommand{\f}[1]{\mathcal{#1}}
\newcommand{\til}[1]{\ensuremath{\widetilde{#1}}}
\newcommand{\m}[1]{\ensuremath{\underline{#1}}}

\newcommand{\proof}{\noindent \textsc{Proof : \ }}

\newcommand{\bo}[1]{\mathbf{#1}}

\newcommand{\pic}{\ensuremath{\mathrm{Pic\ }}}
\newcommand{\D}{\mathfrak{D}}
\newcommand{\const}{\ensuremath{\mathrm{Const}}}
\newcommand{\free}{\ensuremath{\mathrm{Free}}}

\newenvironment{theo} {\begin{list}{$\bullet$}{%
                                \setlength{\leftmargin}{1cm} 
				\setlength{\itemindent}{-0.5cm}
				\setlength{\parsep}{\parskip}
				\setlength{\topsep}{0cm} 
				\setlength{\partopsep}{0cm} 
				\setlength{\itemsep}{0cm} }}
                        {\end{list}}
\newcommand{\hyp}[1]{\ensuremath{#1^\circ\!\textit )}}

\begin{document}
\maketitle

\section*{Introduction}
\addcontentsline{toc}{section}{Introduction}

Let $d,m_1,\ldots,m_r$ be $(r+1)$ positive integers. Denote 
by $V(d;m_1,\ldots,m_r)$
the variety of irreducible (complex) plane curves of degree $d$ having 
exactly $r$ ordinary singularities of multiplicities $m_1,\ldots,m_r$.
In most cases, it is still an open problem to know whether this
variety is empty or not.  

In this paper, we will concentrate on the case where the $r$ singularities
can be taken in a general position. Precisely, let $(P_1,\ldots,P_r)$
be a general $r$-tuple of point in $(\p^2)^r$. Denote by $E$ the
linear system of plane curves of degree $d$ passing through 
the points $P_i$ $(1\leq i\leq r)$ with multiplicity at least $m_i$.
The expected dimension of $E$ is
$\max (-1 ; d(d+3)/2-\sum m_i(m_i+1)/2)$. 

\begin{theorem} 
\label{theorem}
Given a positive integer $m$, there exists an integer $\bo d'(m)$ such
that,
if $m_i\leq m$ for $1\leq i\leq r$ and $d\geq \bo d'(m)$, then~:

The system $E$ has the expected dimension $e$ and, if $e\geq 0$, then 
a general curve in $E$ is irreducible, smooth away from the
$P_i$, and has an ordinary singularity of multiplicity $m_i$ at each
point $P_i$.

As a consequence, $V(d;m_1,\ldots,m_r)$ is not empty.
\end{theorem}

The importance of this result comes from the fact that it is still
valid when the expected dimension
is small (which happens when the number $r$ of points is high) 
; even say, when $e$ is zero. In this case, the curve is 
isolated in $E$, and Bertini's theorem can not be used.

Recent existence results have been proved by Greuel, Lossen and Shustin
in the case of ordinary singularities, (\cite{gls.blowup}, section 3.3) ;
or even for general singularities \cite{gls.plane}. But 
in all these statements the dimension of the system $E$ must be,
at least, quadratic in the 
degree $d$. Notice, however, that the method of \cite{gls.plane}
together with the vanishing result of Alexander-Hirschowitz cited
below (see \cite{al-hi.asymptotic}) would easily give theorem
\ref{theorem} as soon as $e\geq d+1$ (see also section
\ref{sec.highdim} for such considerations).

As for zero dimensional
systems, a previous theorem had been proved by the author 
\cite{mig.mono} for $m_i\leq 3$ and $d\geq 317$.

\medskip

An explicit value for $\bo d'(m)$ has been computed by the author~:
one may take $\bo d'(m)=2((m+2)38)^{2^{m-1}}$. According to a theorem
of Alexander and Hirschowitz \cite{al-hi.asymptotic}, it is already
known that there exists a bound $\bo d(m)$  for the degree, above which
$E$ has the expected dimension. In theorem  \ref{theorem} $\bo
d'(m)$ is slightly greater than $\bo d(m)$ and is expressed in terms
of it. It is now possible to follow the proof of
\cite{al-hi.asymptotic} and give an explicit bound for $\bo d(m)$ 
(let us recall that \cite{al-hi.asymptotic} holds for
any projective variety, since our bound only holds for $\p^2$). With this
approach, it seems that the doubly exponential growth for the explicit
value of $\bo d(m)$ is unavoidable.

However this bound is far from being sharp. In fact, according to a
conjecture of Hirschowitz,
if the $m_i$ are in decreasing order and if $d$ is greater than 
$m_1+m_2+m_3$, then the system $E$ should have the
expected dimension and contain an irreducible and smooth curve 
away from the $P_i$ (except in the
well-known case $(d;m_1,\ldots,m_r)\neq (3n;n,\ldots,n)$ with $r=9$).
Thus, the conjectural bound for $\bo d'(m)$ is
$3m+1$. 

Due to its length, the computation of this explicit value is not
described here. The author places it at the reader's disposal.

\medskip

Theorem \ref{theorem} is also interesting in view of recent results
on the varieties $V(d;m_1,\ldots,m_r)$.
Recall that the first variety of this type, $V=V(d;2,\ldots,2)$
was studied by Severi \cite{sev.anhangf}. He proved that $V$ 
is not empty and smooth if and only if
$r\leq (d-1)(d-2)/2$.  If in addition $r \leq d(d+3)/6$ we also know
that the nodes can be taken
in generic position except in the case $d=6$, $r=9$,
$m_1=\cdots=m_9=2$ (case of an isolated double cubic) (see 
\cite{arb-cor.footnote} and \cite{treg.nodes}).
In 1985, Harris \cite{ha.onseveri} completed this work, proving
that $V(d;2,\ldots,2)$ is always irreducible.

The questions of irreducibility and smoothness 
of general varieties of curves with prescribed singularities
have been treated in many papers. 
Let us mention recent results for general singularities
\cite{sh.equisingular} or for nodal curves on 
general surfaces in $\p^3$ \cite{ch-ci.sevvar}.

However in the case considered here, i.e. plane curves with 
ordinary singularities,
A. Bruno announced that, $V(d;m_1,\ldots,m_r)$ is 
irreducible, smooth and has the expected codimension
assuming that it is not empty and that the singularities can be taken
in generic position (conference in Toledo, September 98). 
This is exactly what is proved in theorem \ref{theorem}.

\medskip

\noindent \emph{\textbf {Strategy of the proof}}

The proof of theorem \ref{theorem} 
is based on a lemma proved by the author in 
\cite{mig.horgeo} (see also \cite{bossini} for a first 
--not differential-- approach of this lemma). This result, which we called
``Geometric Horace Lemma'', is inspired by the Horace method 
of Hirschowitz (see, for example, \cite{al-hi.asymptotic}).
But, while the usual Horace method can only be used to compute
the \emph{dimension} of linear systems
like $E$, the geometrical lemma also yields conclusions about the
\emph{irreducibility} and \emph{smoothness} of the curves in $E$.

The principle of the Geometric Horace Lemma is the following~:
Let us choose an irreducible and smooth plane curve $C$. 
Let us specialize some of the $r$ points on $C$. Denote 
by $y=(Q_1,\ldots,Q_r)$ this special point of $(\p^2)^r$ and by $x$ the
generic point of $(\p^2)^r$. Two linear systems may be considered~:
$E_x=E$ when the points are in generic position and  $E_y$ when
they are in special position.
The specialization from $x$ to $y$ is done in such a way, that 
$C$ is a base component of the system $E_y$.
Thus a curve in $E_y$ is the union of $C$ and of a \emph{residual}
curve.

Under some assumptions, detailed in \ref{horgeo}, if the generic
residual curve is geometrically irreducible, smooth, and has ordinary
singularities, then the general curve in $E_x$ also satisfies these
properties. 

An important point must be mentioned~:
if we do not specialize enough points on $C$, then $C$ is not
 a base component of $E_y$ and the method fails. But, if we 
specialize too
many points, then the dimension of the linear system grows~:
$\dim E_y > \dim E_x$. This phenomenon is controlled 
with the help of differential conditions. It means that we
have to consider some sub-systems of curves bound to pass
through infinitely near points.

\medskip

Here is the main point of the proof~: by specializing too many points
on the curve $C$, it is possible to make the dimension of $E$
grow considerably~; i.e. grow as high as the degree $d$. Then, 
assuming that some vanishing property holds true,
the residual system is base point free, 
and Bertini's theorem can be used. As a consequence, a
general residual curve is smooth, irreducible, and 
the intersection variety described above is irreducible.

\medskip

To make all this strategy work, we still have to check the vanishing
property referred to above. Roughly speaking, it means that the
residual system has the expected dimension. 
To prove this, we make use of the following 
vanishing result of Alexander and Hirschowitz \cite{al-hi.asymptotic}~:

Given an integer $m$, there exists an integer $\bo a(m)$, and for
$a\geq \bo a(m)$, there exists another bound $d_0(a,m)$ such that, 
if $C$ is the generic curve of degree $a$, if $d\geq d_0(a,m)$ and 
if the points $P_i$ are either generic in $\p^2$ or generic on $C$ 
(not too many of them) then the system $E$ has the expected dimension.

In view of this result, the last choices are made~: As for the
curve $C$, we choose the generic curve of degree $\bo a(m)$ ;
and we only consider systems of curves of degree higher than
$\bo d(m)=d_0(\bo a(m),m)$  (in fact, the final value $\bo d'(m)$ is
greater than $\bo d(m)$, as appears in theorem
\ref{theorem2}).

\medskip

\noindent \emph{\textbf {Contents}}

The article is organized as follows~:
In the first part, notations and definitions are set. In particular, we
describe the universal variety which parameterizes the curves we are
studying. 

In the second and third sections, we restate respectively 
the Geometric Horace Lemma, and
the vanishing theorem of Alexander and Hirschowitz. These are the two
main tools in the proof of theorem \ref{theorem}.

The fourth section is devoted to the study of linear systems of ``high''
dimension, (precisely, a dimension greater than the degree $d$).
In particular, when the $r$ points are in a good position, 
so that the vanishing lemma can be used, we show 
that theorem \ref{theorem} is true for these systems.

In the last section , theorem \ref{theorem} is proved.

\section{Curves on rational surfaces}
\label{sec.prelim}

In the introduction, the situation has been described on the plane.
Actually, most of the proofs will be done on the plane blown-up along the
$r$ points $P_1,\ldots,P_r$. 

In this section, we shall describe the family of rational surfaces obtained
by blowing up a family of $r$ disjoint sections in $\p^2$, and the
families of curves on these surfaces. We shall also set up most of the
notations used in the article.

\subsection{Families of rational surfaces and relative divisors}

Let $r$ be a positive integer, and $X\subset (\p^2)^r$ be the open
subset of $r$-tuples of distinct points. The morphism
$\p^2_X=\p^2\times X\fl X$ is naturally endowed with $r$ sections~:
$$
\begin{array}{cccc}
\gamma_i : & X &                \fl &           \p^2\times X\\
&(P_1,\ldots,P_r)&      \fl&            (P_i,(P_1,\ldots,P_r))
\end{array}
$$
Let $\Gamma_i$ be the image of $\gamma_i$ ; $\Gamma=\cup_{i=1}^r
\Gamma_i$ is a nonsingular variety of $\p^2_X$. Blowing up $\p^2_X$
along $\Gamma$ produces a family of rational surfaces, parameterized by
$X$~: $S_X\morph{b}\p^2_X\fl X$. Let $\pi$ denote 
the composed morphism $S_X\fl X$. 
At any point $x=(P_1,\ldots,P_r)$ of $X$, the fiber of $\pi$ will be
denoted by $S_x$. This surface $S_x$ is simply the projective plane blown up
along the $r$ points $P_1,\ldots,P_r$. 

\medskip

Let us keep in mind that a relative effective Cartier divisor of $S_X$
on $X$ is simply an ideal sheaf $\f I_D$ on $S_X$, locally principal, and 
not a zero-divisor in  any fiber of $\pi$ (see \cite{gro.pic}). These
ideal sheaves are flat on $X$.

\medskip

\noindent \textbf{Examples}~: \hyp 1 Consider a line $L$ on $\p^2$,
$L\times X\subset \p^2\times X$ the trivial family of lines 
above $X$, and $H_X=b^{-1}(L\times X)$ the total transform of $L\times
X$ in $S_X$. The ideal sheaf $\f I_{H_X}$ is a relative effective
Cartier divisor of $S_X$ on $X$.

\hyp 2 Consider now $E_{i,X}$, the exceptional divisor obtained by 
blowing up the irreducible smooth variety $\Gamma_i$ ; the ideal sheaf
$\f I_{E_{i,X}}$ is also a relative effective
Cartier divisor of $S_X$ on $X$.

\subsection{Intersection pairing and linear systems}

For any $x\in X$, the Picard group of the surface $S_x$ is endowed
with the usual base~: 
$[H_x],[-E_{1,x}],\ldots,[-E_{r,x}]$. The relative Cartier divisors
being flat on $X$, these bases satisfy the following property~:
\begin{proposition}
Let $D$ be a relative Cartier divisor of $S_X\fl X$. For any point
$x\in X$, let $[D_x]\in \Z^{r+1}$ be the class of the sheaf $\f
O_{S_x}(D)$ expressed in the base defined above ;
then $[D_x]$ does not depend on $x$.
\end{proposition}

The canonical divisor of $\pic S_x$ is $\m \omega_x=(-3;(-1)^r)$
(the notation $(-1)^r$ means that the integer $(-1)$ is repeated $r$ 
times ; this convention will be kept in the sequel). The intersection
pairing of $\pic S_x$ is as follows~:
$(d;m_1,\ldots,m_r).(c;n_1,\ldots,n_r)=dc-m_1n_1-\ldots -m_rn_r$.

Let $\m d=(d;m_1,\ldots,m_r)\in \pic S_x$. The 
sheaf $\f O(dH_x-\sum_{i=1}^r m_iE_{i,x})$ will be denoted by 
$\f O(\m d)$, and the complete linear system of $\f O(\m d)$ 
(i.e. the projective space $\p(H^0(S_x,\f O(\m d)))$ will be 
denoted by $\f L_x(\m d)$. 

The Riemann-Roch theorem for surfaces allows us to compute the
Euler-Poincar\'e characteristic of $\f O(\m d)$~:
$\chi(\m d)=\frac{(d+1)(d+2)}{2}-
\sum_{i=1}^r\frac{m_i(m_i+1)}{2}.$
One may also compute the arithmetical genus of the eventual sections of
$\f O(\m d)$~:
$g(\m d)=\frac{(d-1)(d-2)}{2}-
\sum_{i=1}^r\frac{m_i(m_i-1)}{2}.$

Suppose that $d\geq -2$. By Serre Duality Theorem,
$h^2(S_x,\f O(\m d))=0$. The \emph{expected dimension}
for $\f L_x(\m d)$ is then~:
$\max(\chi(\m d)-1,-1)$. (An empty system is supposed to have dimension
$-1$). A system $\f L_x(\m d)$ will be said to be \emph{regular}
if it has the expected dimension. More generally, a sheaf $\f F$ on
$S_x$ will be said to be regular if $h^0(S_x,\f F)=
\max(\chi(\f F),0)$ and $h^1(S_x,\f F)=
\max(-\chi(\f F),0)$.

\subsection{Universal family of divisors}
\label{family}

Let $\m d=(d;m_1,\ldots,m_r)$ be an $r$-tuple of integers, and 
let us define  
$m'_i=\max(m_i,0)$ and $\m d'=(d;m'_1,\ldots,m'_r)$. 
Let $\Gamma_i$ be
the image of the $i^{th}$ natural section $\gamma_i$ defined above,
and let $Z\subset \p^2_X$ be the scheme defined by the ideal
$\f I_{\Gamma_1}^{m'_1}\ldots \f I_{\Gamma_r}^{m'_r}$.
The scheme $Z$ is a flat family above $X$ whose fibers $Z_x$ are
unions of $r$ fat points of multiplicities $m'_1,\ldots,m'_r$
(see section \ref{sec.vanish} for a definition of fat points).

Consider $x=(P_1,\ldots,P_r)\in X$. The linear system $|\f
I_{Z_x}(d)|$ is the system of plane curves of degree $d$ passing
through each point $P_i$ with multiplicity at least $m'_i$.
One can easily see that this system is isomorphic to $\f L_x(\m d)$. 

\medskip

Consider now the linear system $\p(H^0(\p^2,\f O(d)))\isom
\p^{d(d+3)/2}$ of plane curves of degree $d$. One may define
(here, only under a ``set-theoretic'' point of view, but it is endowed
with a natural scheme structure) a subscheme $F$ of
$\p^{d(d+3)/2}\times X$ in the following way~:
\smallskip

$\begin{array}{rcll}
F& = & \{(D,(P_1,\ldots,P_r))\  | & \mbox{$D$ passes through each point} \\
	& &	&\mbox{$P_i$ with multiplicity at least $m'_i$}\}\\
\end{array}$

\smallskip
 
This scheme $F$ parameterizes a canonical family of curves $\f
D'\subset \p^2\times F$~: given $x=(D,(P_1,\ldots,P_r))$, the fiber
$\f D'_x$ simply is the curve $D$. As above, it is possible to blow-up
the variety $\p^2_F=\p^2\times F$ along the $r$ disjoint natural 
sections. Let $b_F:S_F\fl \p^2_F$ be this blowing-up. By assumption,
the divisor $\f D=b_F^{-1}(\f D')-\Sigma_{i=1}^r m_iE_i$ is
effective~: it is a relative effective Cartier divisor of the family
$S_F\fl F$. Moreover, $\f D$ is a universal divisor~:
\begin{proposition}
Let $\m d\in \Z^{r+1}$. The functor $\bo F_{\m d}$ from
$\mathbf{Schemes/X}$ to $\mathbf{Sets}$ such that
$\bo F_{\m d}(Y)=\{D\subset S_Y, \mbox{relative effective divisor of
class $\m d$ on $Y$}\}$ is represented by the couple $(F,\f D)$.
\end{proposition}
(This proposition is detailed in \cite{mig.these} ; see also
\cite{gro.pic} and \cite{no.planecur}).

\medskip 

Let $p:F\fl X$ be the natural projection from $F\subset \p^{d(d+3)/2}
\times X$ to $X$. The fiber of $p$ over a point $x\in X$ is nothing 
but the linear system $\f L_x(\m d)$. Let $x'$ be the generic point of 
this fiber. By definition of the universal divisor $\f D$, the curve
$\f D_{x'}$ is the generic curve of $\f L_x(\m d)$.
It will be denoted by $\D_x(\m d)$.

Suppose now that a point $y\in X$ is a specialization of $x$. 
We will say that the curve $\D_x(\m d)$ \emph{specializes to the
curve} $\D_y(\m d)$ if the generic point of $p^{-1}(y)=\f L_y(\m d)$
is a specialization of $x'$. This notion is of special importance 
if one expects to find properties of $\D_x(\m d)$ 
from those of $\D_y(\m d)$. In particular, if $\D_y(\m d)$
is geometrically irreducible, smooth, or has ordinary singularities, 
then the same holds for $\D_x(\m d)$. 

If the dimension of $\f L_y(\m d) $ equals the dimension 
of $\f L_x(\m d)$, one can easily see that  $\D_x(\m d)$ always
specializes to $\D_y(\m d)$. If the dimension grows after the
specialization, extra conditions are needed. In fact, it is sufficient
to prove that the strata of the cohomological stratification 
(associated to the sheaf $\f O_{S_X}(\m d)$) have sufficiently 
big codimension. This can be done with the help of differential
methods (see \cite{mig.horgeo}, and the lemma \ref{horgeo}).

\section{The Geometric Horace Lemma}
\label{sec.horgeo}

In this section, the Geometric Horace Lemma is restated and
commented on. This lemma was
proved by the author in \cite{mig.horgeo}. 

\medskip

Let us first give some notations and conventions~:
Let $y=(Q_1,\ldots,Q_r)$ be a point of
$X$ (the notation $(Q_1,\ldots,Q_r)$ is slightly
incorrect, since $y$ is generally not a closed point).
Let $G$ be a closed integral subscheme of $\p^2$. We will say
that the $a$ points $Q_1,\ldots, Q_a$ $(0\leq a\leq r)$ are 
\emph{generic and independent} on $G$ if $y$ is the generic point of a
subvariety $Y\subset X$ such that $Y=G^a\times V$,
where $V$ is an irreducible subscheme of $(\p^2)^{r-a}$.

Suppose now that the $i$-th point $Q_i$ 
is a nonsingular point of a plane curve $C$. On the rational
surface $S_y$, the intersection point of 
the exceptional divisor $E_i$ and the strict transform $\til C$ will
be denoted by $Q_i^C$ (and its ideal sheaf, $\f I_{Q_i^C}$).
Let $\f O(\m d=(d;m_1,\ldots,m_r))$ be an invertible sheaf on $S_x$.
Global sections of the sheaf 
$\f I_{Q_i^C}(\m d)$  can be seen as plane curves of
degree $d$ having multiplicity at least $m_j$ at each point $Q_j$
and, if the multiplicity at $Q_i$ is exactly $m_i$
having a branch tangent to $C$ at this point.

\begin{lemma}
\label{horgeo}
Let $\m d=(d;m_1,\ldots,m_r)$ be an $r$-tuple of positive integers, 
and $x=(P_1,\ldots,P_r)$, $y=(Q_1,\ldots,Q_r)$ be two points of $X$ 
such that $x$ specializes to $y$. 
Let $C$ be a plane curve, and $\til C:=\til C_y$  its strict transform on
$S_y$. Assume that $\til C$ is geometrically irreducible and smooth,
of class $\m c\in \pic S_y$ and of genus $g(\m c)$. 

\smallskip

\noindent \textbf{Dimension and specialization} Suppose that~:

\begin{theo}
\item[$\hyp 1$]
$-\al:=\chi(\f O_{\til C}(\m d))=\m d.\m c+1-g(\m c)\leq 0$\ \ ;
\item[$2^\circ )$] At the point $y$, $g(\m c)$ points
are generic and independent on $C$.
\item[$3^\circ )$] If $\al \geq 1$, there exist 
$\al +1$ integers $i_1,\ldots,i_{\al +1}$ such that~: \\
$P_{i_1},\ldots,P_{i_{\al +1}}$ are generic and independent in the
plane, and \\
$Q_{i_1},\ldots,Q_{i_{\al +1}}$ are generic and independent on $C$.
\item[$4^\circ )$] If $\al =0$, $H^0(S,\f O(\m d-\m c))$ has the expected
dimension~: $\chi(\m d-\m c)=\chi(\m d)$.\\
If $\al \geq 1$, 
$H^0(S,\f I_{Q_{i_1}^C\cup\cdots\cup Q_{i_{\al +1}}^C}(\m d-\m c))$ 
has the expected dimension~: $\chi(\m d-\m c)-\al -1 = 
\chi(\m d)-1$,
\end{theo}
\noindent Then $\f L_x(\m d)$ is regular, and, if $\chi(\m d)>0$, 
$\D_x(\m d)$ specializes to $\D_y(\m d)$.

\noindent \textbf{Irreducibility}
Suppose, moreover, that $\chi(\m d)>0$ and~:
\begin{theo}
\item [$\hyp 5$] the system $\f L_x(\m c)$ is empty.
\item [$\hyp 6$] $\D_y(\m d-\m c)$ is geometrically irreducible. 
\end{theo}
\noindent Then $\D_x(\m d)$ is geometrically irreducible.

\noindent \textbf{Smoothness}
Suppose finally that~:
\begin{theo}
\item [$\hyp 7$] If $\al =0$, $y$ is normal and the closure of $x$~;
\item [$\hyp 8$] $\D_y(\m d-\m c)$ meets $\til C$ transversally~;
\item [$\hyp 9$] $\D_y(\m d-\m c)\cap \til C$ is irreducible (not: 
geometrically irreducible)~;
\item [$\hyp {10}$] $\D_y(\m d-\m c)$ is a smooth curve,
\end{theo}
then $\D_x(\m d)$ is smooth.
\end{lemma}

The system $\f L_y(\m d-\m c)$, and the curve $\D_y(\m d-\m c)$ are
respectively called \emph{residual system} and  \emph{residual curve}.
In fact, $\D_y(\m d)$ is the union of $\til C$ and
$\D_y(\m d-\m c)$. The curve $\til C$ being well-known, this lemma
allows one to get information on $\D_x(\m d)$,
from the properties of $\D_y(\m d-\m c)$ and the relation between
$\D_y(\m d-\m c)$ and $\til C$.

\medskip

Condition \hyp 9 is certainly best described in an example~:
assume that $\dim \f L_y(\m d)=0$, and denote by $Y$ the adherence of $y$
in $X$. There is only one curve in the system 
$\f L_z(\m d)$ for every closed point $z$ in an open subset of $Y$. 
Then the intersection of $\til C$ and the residual curve in $\f L_z(\m
d-\m c)$ makes, as $z$ varies, a covering of degree $(\m d-\m c).\m c$ over 
the open subset of $Y$. The condition \hyp 9 means that this
``intersection variety'' is irreducible ; or in other words, 
that the monodromy group of this covering acts transitively
on the intersection points of $\til C$ and
the residual curve $\D_y(\m d-\m c)$.

\medskip

A case of special interest is the case where the number $\al$ of
\hyp 1 is positive. Roughly speaking, this situation arises when one
specializes too many points on the curve $C$. Let us suppose that
$\chi(\m d)>0$. Considering the exact sequence
\begin{eqnarray}
\label{exact}
0\fl \f O_{S_y}(\m d-\m c)\fl \f O_{S_y}(\m d)\fl \f O_{\til C}(\m
d)\fl 0,
\end{eqnarray}
we find that $\chi(\m d-\m c)=\chi(\m d)+\al $. Since
$\f I_{Q_{i_1}^C\cup\cdots\cup Q_{i_{\al +1}}^C}(\m d-\m c) $ 
is regular, $\f L_y(\m d-\m c)$ also is regular
and $\dim \f L_y(\m d)=\dim \f L_y(\m d-\m c)=\dim \f L_x(\m d)+\al $.
Thus, the dimension has grown by $\al$.

\begin{remark}
\label{ordinary}
\textnormal{
As regards to the ordinary singularities of the plane projection of 
$\D_x(\m d)$, the Geometric Horace Lemma
yields no conclusion. But, if $\D_y(\m d)=\D_y(\m d-\m c)\cup \til C$
meets the exceptional divisor $E_i$ in $m_i$ distinct points,
$\D_x(\m d)$ also possesses this property. Hence, its projection
has only ordinary singularities of the expected multiplicity.
}

\textnormal{
This is in particular the case if $C$ has ordinary singularities and
$\f L_y(\m d-\m c)$ is base point free
(Bertini's theorem applied to $E_i$ and the restricted system
$\f L_y(\m d-\m c)_{|E_i}$).
}
\end{remark}

\section{An asymptotic vanishing theorem}
\label{sec.vanish}

In order to use the Geometric Horace Lemma, we have to check that some
linear system is regular (condition \hyp 4 of \ref{horgeo}). 
This will be done with the
help of an asymptotic vanishing theorem of Alexander and Hirschowitz.
In this section, we restate this result and the adequate definitions.
As a corollary, we write down precisely the vanishing lemma used in
the proof of theorem~\ref{theorem}. 

Here, opposed to the other sections, 
all the work is done on the plane, without
blowing it up. It is a natural choice when dealing
with the \emph{dimension} of a linear system, without consideration
of the smoothness of its sections.

\medskip

Let us first recall some definitions~:
As usual, a \emph{fat point} of support
$P\in \p^2$ is a subscheme $P^m$ of $\p^2$ defined by the ideal $\f I_P^m$ ;
the integer $m$ is called the \emph{multiplicity } of $P$. 
If $Z$ is a zero dimensional subscheme of $\p^2$, the \emph{degree}
of $Z$, denoted by $\deg Z$, is the length of the ring $\f O_Z$.
As an example, $\deg P^m=m(m+1)/2$.
\begin{definition} Let $P\in \p^2$ be a nonsingular point of a plane
curve $C$, and $i,m$ be two integers such that $0\leq i\leq m-1$.

The \emph{$i$-th residue point} 
supported by $P$, of multiplicity $m$, with respect 
to $C$, is the scheme defined by the ideal
$\f I_P^{m-1}\cap(\f I_C^i+\f I_P^m)$. 
It is denoted by 
$D_{C}^i(P^m)$ or $D^i(P^m)$ if no confusion can arise. 
A residue of type $D^{m-1}(P^m)$ is called a \emph{simple residue}
(\cite{al-hi.asymptotic}, 2.2.).
\end{definition}

Condition \hyp 4 of lemma \ref{horgeo} is easily expressed
with the help of residue points~: For example, 
if $P_1,\ldots,P_r$ are points of $\p^2$, such that 
$P_1$ is a nonsingular point of a curve $C$, and if
$\m d=(d;m_1,\ldots,m_r)$ is an $(r+1)$-tuple of positive integers, 
we consider the zero dimensional scheme 
$Z=D^1(P_1^{m_1+1})\cup P_{2}^{m_{2}} \cup\cdots\cup P_r^{m_r}$.
Then $|\f I_{P_1^C}(\m d)|$ is isomorphic to $|\f I_Z(d)|$. 
What we still have to show is that $h^1(\p^2,\f I_Z(d))$ or
$h^0(\p^2,\f I_Z(d))$ is zero. 

The vanishing result of Alexander and Hirschowitz deals with some
special types of systems $|\f I_Z(d)|$ defined below~:

\begin{definitions}[\cite{al-hi.asymptotic}, 3.1]
\label{defcand} 
Let $d,m$ and $a$ be three positive integers ; denote by $C$ the generic
plane curve of degree $a$.

An \emph{$(m,a)$-configuration} is a zero dimensional scheme
$Z=\const(Z)\cup\free(Z)$, where~: $\free(Z)$ is the \emph{free} part
of $Z$ ; it is a union of fat points,
of generic and independent support in $\p^2$.  $\const(Z)$ 
is the \emph{constrained} part of $Z$ ; it is a union 
of fat points or simple
residue points, of generic and independent support in $C$.
All these points are supposed to have a multiplicity less than or
equal to $m$.

A \emph{$(d,m,a)$-candidate}
is an $(m,a)$-configuration $Z$ such that
$\chi(\f I_Z(d))\leq 0$ and
$h^0(C,\f O_{C}(d))\geq \deg(Z\cap C)$.		

A $(d,m,a)$-candidate such that $H^0(\p^2,\f I_Z(d))=0$ is said to be
\emph{winning}. 
\end{definitions}

\begin{remark}
\label{remark}
\textnormal{
If $Z$ is an $(m,a)$-configuration such that $\chi(\f I_Z(d))>0$ and
$h^0(C,\f O_{C}(d))\geq
\deg(Z\cap C)$, one also says that
$Z$ is a $(d,m,a)$-candidate. But in this case, $Z$ is winning
means that $h^1$ vanishes, whereas $h^0$ is positive.
}
\end{remark}

The crucial vanishing result is the following~:

\begin{proposition}[\cite{al-hi.asymptotic}, 7.1]
\label{theoalhi}
Given a positive integer $m$, 
there exists an integer $\bo{a}(m)$ and, for each $a\geq \bo{a}(m)$,
an integer $d_0(a,m)$ such that :
if $a\geq \bo{a}(m)$ and $d\geq d_0(a,m)$ then any 
$(d,m,a)$- candidate is winning.   
\end{proposition}

Unfortunately, this vanishing result 
involves \emph{simple residues}, of
type $D^{m-1}$, whereas the needed condition involves residues of
type $D^1$. With our ``asymptotical'' point of view, 
this is essentially a technical problem. 
But to solve it, a little more Horace method is needed~:

Let $Z$ be a closed subscheme of $\p^2$, and $C$ be an irreducible and
reduced plane curve. The \emph{trace} of $Z$, denoted $Z\cap C$, is the
scheme defined by the ideal $\f I_Z+\f I_C$. The residue of
$Z$ with respect to $C$,  denoted $Z'$, is the
scheme defined by the conductor ideal $(\f I_Z:\f I_C)$. 

\begin{proposition}[\cite{al-hi.asymptotic}, 2.3]
\label{hordiff}
Let $C$ be an irreducible and smooth plane curve of degree $a$ and
genus $g=(a-1)(a-2)/2$, and
$d$ be an integer greater than $a$. Let $Z=Z_0\cup P_1^{m_1}\cup\cdots
\cup P_\beta ^{m_\beta }$ be a zero dimensional subscheme of $\p^2$ such that 
$P_1,\ldots,P_\beta $ are generic and independent points of $\p^2$. 
Denote also by $Q_1,\ldots,Q_\beta $, $\beta $ generic and independent points of
$C$. Suppose that $\chi(\f I_Z(d))\leq 0$. If~:

$i)$ $H^1(\p^2,\f I_{Q_1^{m_1}\cup\cdots\cup Q_\beta ^{m_\beta
}}(d-a))=0$ ;

$ii)$ $\beta  = da+1-g-\deg(Z\cap C)$ ;

$iii)$ $\f I_{(Z_0\cap C)\cup Q_1 \cup\cdots\cup Q_\beta }(d)$ is a regular
invertible sheaf of $C$ ;

$iv)$ $H^0(\p^2,\f I_{Z'_0\cup D^{m_1-1}(Q_1^{m_1})\cup\cdots\cup 
D^{m_\beta -1}(Q_{\beta}^{m_\beta })}(d-a))=0$ ; 

\noindent then $H^0(\p^2,\f I_Z(d))=0$, as expected.
\end{proposition}

\begin{corollary}
\label{vanish}
Let $m,a$ be two positive integers, and $C$ be the generic plane curve
of degree $a$. Denote by $x=(O_1,\ldots,O_t,P_1,\ldots,P_r)$ the generic point of
$C^t\times (\p^2)^r$, and consider an integer $\al$ such that 
$0\leq \al \leq t$. Let $\m d=(d;n_1,\ldots,n_t,m_1,\ldots,m_r)\in \pic
S_x$ such that $n_i\leq m-1,\  m_i\leq m$, and 
suppose that $\chi(\m d)-\al= 0$. If

$i)$ $a\geq \max(\bo a(m)\ ,\ 4m)$ ;

$ii)$ $d\geq \max(d_0(a,m)+a\ ,\ 2am)$ ;

$iii)$ $da+1-g-n_1-\ldots-n_t-\al \geq 0$ ;

\noindent then $\f I_{O_1^C\cup\ldots\cup O_{\al}^C}(\m d)$ 
is regular.
\end{corollary}

\proof Let $Y_0=D^1(O_1^{n_1+1})\cup\cdots\cup
D^1(O_\al^{n_\al+1})\cup O_{\al+1}^{n_{\al+1}}
\cup\cdots\cup O_t^{n_t}$,
and $Y=Y_0\cup P_1^{m_1}\cup\cdots\cup P_r^{m_r}$. Clearly,
$\chi(\f I_Y(d))=\chi(\m d)-\al=0$, therefore 
$\f I_{O_1^C\cup\ldots\cup O_{\al}^C}(\m d)$ is regular if and only if
$H^0(\p^2,\f I_Y(d))=0$.

\medskip

\noindent $\bullet$ Let us first prove that
there exists a non negative integer $s$ such that~:

\vspace{-0.5cm}

\begin{eqnarray}
\label{ajust}
\beta&:=& da+1-g-\al-\sum_{i=1}^tn_i-\sum_{j=1}^sm_j\in [0;m-1]\\
\label{derivable}
s+\beta&\leq r
\end{eqnarray}

\vspace{-0.1cm}

Since $0<m_i\leq m$, it will be enough to show that~:
$$
\begin{array}{ll}
\ &\!\! da+1-g-\al-\sum_{i=1}^t n_i-\sum_{j=1}^{r-m+1}m_j\leq m-1\\
\Longleftarrow_{(m_j\leq m)}&\!\! 
	da+1-g-\al-\sum_{i=1}^t n_i-\sum_{j=1}^{r} m_j + (m-1)m\leq m-1\\
\Longleftarrow_{(m_j,n_i\leq m)}&\!\! da+1-g-\al\\
	&\mbox{\qquad}-\frac{2}{m+1}\left(\sum \frac{n_i(n_i+1)}2 	
	+\sum\frac{m_i(m_i+1)}2\right) +(m-1)^2 \leq 0\\
\equi &\!\! da+1-g-\al-\frac{2}{m+1}\left(\frac{(d+1)(d+2)}2-\chi(\m
	d)\right) + (m-1)^2 \leq	0\\ 
\equi_{(\chi(\m d)= \al)} &\!\! da+1-g-\frac{(d+1)(d+2)}{(m+1)}-\al
	\left(1-\frac{2}{m+1}\right) + (m-1)^2 \leq 0\\
\Longleftarrow_{(m>0)} &\!\! da+1-\frac{(a-1)(a-2)}{2}
	-\frac{(d+1)(d+2)}{(m+1)} + (m-1)^2 \leq 0\\
\Longleftarrow_{\scriptscriptstyle{(\!d+1\geq a(\!m+1\!)\!)}} &\!\!
1-\frac{(a-1)(a-2)}{2} 
	-2a + (m-1)^2 \leq 0\\
\Longleftarrow_{(a\geq 4 m)} &\!\! 1-
	\frac{(4 m-1)(4 m-2)}{2}
	-8 m + (m-1)^2 \leq 0\\
\equi &\!\! -7m^2-4m+1\leq 0, \mbox{\ \ which is true.}
\end{array}
$$

\medskip

\noindent $\bullet$
Consider $Q_1,\ldots,Q_{s+\beta}$, $(s+\beta)$ generic
and independent points on $C$. Denote by $Z_0$ and $Z$ the schemes
\begin{eqnarray*}
Z_0&:=& Y_0\cup Q_1^{m_1}\cup\cdots\cup Q_s^{m_s}\cup
	P_{s+\beta+1}^{m_{s+\beta+1}}\cup \cdots\cup P_r^{m_r}\\
Z &:=& Z_0\cup P_{s+1}^{m_{s+1}}\cup \cdots\cup P_{s+\beta}^{m_{s+\beta}}.
\end{eqnarray*}
By the Semicontinuity Theorem, if $H^0(\p^2,\f I_Z(d))=0$, then 
$H^0(\p^2,\f I_Y(d))=0$ as expected.

To prove that $H^0(\f I_Z(d))$ is equal to zero, we make use of proposition
\ref{hordiff}.
Let us specialize the points $P_{s+1},\ldots,P_{s+\beta}$
to the points $Q_{s+1},\ldots,Q_{s+\beta}$.
The following relation,which bound the number of generic points on
$C$, will be useful~:
\begin{equation}
\label{bigts}
t+s+\beta \geq 2a^2-\frac{a^2}{2m}
\end{equation}
This inequality comes from (\ref{ajust}), which yields
$\sum_{i=1}^t n_i + \sum_{j=1}^{s}m_j +\beta \geq da+1-g$.
Since $n_i,m_j\leq m$ and $d\geq 2am$  one gets 
$m(t+s+\beta)\geq 2a^2m-a^2/2 +3a/2$.

Let us check conditions $i$ to $iv$ of \ref{hordiff}~:

$i)$ $H^1(\p^2,\f I_{Q_{s+1}^{m_{s+1}}\cup\cdots\cup
Q_{s+\beta}^{m_{s+\beta}}}(d-a))=0$ by the lemma \ref{xu} below. 

$ii)$ By definition of $Z$, $\deg Z\cap C=
\al+\sum_{i=1}^t n_i+\sum_{j=1}^s m_j$. So that
$\beta=da+1-g-\deg Z\cap C$.

$iii)$ The divisor of $C$ defined by the ideal
$\f J=\f I_{(Z_0\cap C)\cup Q_{s+1} \cup\cdots\cup Q_{s+\beta}}$ 
is supported on the $t+s+\beta$ points
$O_1\ldots,O_t,Q_1,\ldots,Q_{s+\beta}$ which are generic and independent
on $C$. Hence, if $t+s+\beta\geq  g$, $\f J(d)$ is a nonspecial
invertible sheaf. But, $t+s+\beta \geq 2a^2-a^2/(2m)$
(\ref{bigts}) and then $t+s+\beta \geq a^2/2\geq g$.

$iv)$ Let $T=Z'_0\cup D^{m_{s+1}-1}(Q_{s+1}^{m_{s+1}}) \cup\cdots\cup
D^{m_{s+\beta}-1}(Q_{s+\beta}^{m_{s+\beta}})$. 
Since $Z'_0$ is a union of fat points
(the $D^1$ have disappeared), $T$ is an $(m,a)$-configuration.
Let us prove that it is a $(d-a,m,a)$-candidate.

From the definition of $\beta$, one easily sees that
$\chi(\f I_T(d-a))=\chi(\f I_Z(d))=0$. The second condition is~:
$h^0(C,\f O_C(d-a))-\deg (T\cap C) \geq 0$.
If $s+1\leq j\leq s+\beta$, then 
$\deg (D^{m_j-1}(Q_j^{m_j}))$ equals $m_j$ if $m_j> 1$ and $0$
if $m_j=1$. So the inequality can be checked as follows~:
$$
\begin{array}{ll}
& h^0(C,\f O_C(d-a))-\deg (T\cap C)\\
\geq & (d-a)a+1-g-\sum_{i=1}^t (n_i\!-\!1)-\!
\sum_{j=1}^s(m_j\!-\!1)\!-\sum_{j=s+1}^{s+\beta}m_j \\
=_{(\ref{ajust})}& -a^2+t+s+\al+\beta - \sum_{j=s+1}^{s+\beta}m_j \\
\geq_{(m_j\leq m)}& -a^2+t+s+\al+\beta -m(m-1)\\
\geq_{\mathrm{(\ref{bigts})}}& -a^2+(2a^2-a^2/2m) -m(m-1)\\
\geq_{(a\geq 4 m)}& 
	15m^2 -7m \geq 0.
\end{array}
$$
Thus $T$ is a $(d-a,m,a)$-candidate. By assumption, 
$a\geq \bo a(m)$, and $d-a\geq d_0(a,m)$ ; hence, by Proposition 
\ref{theoalhi}, $T$ is winning and $H^0(\f I_T(d-a))=0$.

\smallskip

It is now allowed to apply proposition \ref{hordiff} ; it gives
$H^0(\f I_Z(d))=0$. \sq

\begin{lemma}
\label{xu}
Under the assumptions of corollary \ref{vanish},
consider $Q_1,\ldots,Q_\beta$, $\beta $ generic points of
$C$, and $m_1,\ldots,m_\beta$, $\beta $ integers bounded by $m$.
If $\beta \leq (m-1)$, then $H^1(\p^2,\f I_{Q_1^{m_1}\cup\cdots\cup
Q_\beta^{m_\beta}}(d-a))=0$. 
\end{lemma}
\proof
The only case we need to consider is $\beta=m-1$ and
$m_1=\ldots=m_\beta=m$. By Xu's theorem 
(\cite{xu.ample}, theorem 3),  
$H^1(\p^2,\f I_{Q_1^{m}\cup\cdots\cup Q_\beta^{m}}(d-a))=0$,
as soon as 
$d-a\geq 3m$ and $(d-a+3)^2>(10/9)\sum^{\beta}_{i=1}
(m_i+1)^2=(10/9)(m-1)(m+1)^2$.
By assumption, $d-a\geq 2am-a\geq am$ and $a\geq \sqrt 2 m$, hence
$(d-a+3)^2\geq a^2m^2 \geq 2m^4 \geq (10/9)(m-1)(m+1)^2$
for any positive $m$. If $a\geq 3$ then 
$d-a \geq 3m$. If $a=2$ then $m$ necessary equals $1$, and the lemma 
is clearly true. 
\sq

\section{Systems of high dimension}
\label{sec.highdim}

Given a system $\f L(\m d)$ of a sufficiently high dimension,
Theorem \ref{theorem} may be proved with the help of Bertini's
theorem. The main point consists in showing that $\f L(\m d)$ is base
point free ; this is done here with the vanishing theorem of the preceding
section and a kind of Castelnuovo-Mumford's argument.

\begin{proposition}
\label{highdim}
Let $m$ be a positive integer, $a\geq \bo a(m)$ (see prop.
\ref{theoalhi}) and $C$ the generic curve of degree
$a$. Let $x=(P_1,\ldots,P_r)$ be the generic point of $C^s\times
(\p^2)^{r-s}$, $(0\leq s\leq r)$ and $\m d=(d;m_1,\ldots,m_r)\in \pic
S_x$ be such that $m_i\leq m$ ($1\leq i\leq r$). 
The class of $\til C$ in $\pic S_x$ is denoted by $\m c$, and its genus $g$.

Suppose that $d\geq d_0(a,m)+1$ (see \ref{theoalhi}), $\chi(\m d)\geq d+1$,
and $\m d.\m c+1-g \geq a$. Then,
$\f L_x(\m d)$ is base point free,
$\D_x(\m d)$ is geometrically irreducible and smooth,
$\D_x(\m d)$ meets $\til C$ transversally, and $\D_x(\m d)\cap
\til C$ is irreducible.
\end{proposition}
\proof Let us first prove that $\f L_x(\m d)$ is base
point free~: the characteristic $\chi(\m d)$ is greater than $1$, so
we just have to show that, given a point $Q\in S_x$, 
$h^1(S_x, \f I_Q(\m d))=0$.
Whatever the position of $Q$ is (even on an exceptional divisor),
there exists a line $L$ on $\p^2$ such that $Q$ belongs to the strict
transform $\til L$ of $L$. Let $\m l\in \pic S_x$ be the class of
$\til L$ ; one may write $\m l=(1;\eps_1,\ldots,\eps_r)$, where
$\eps_i=1$ if $P_i\in L$ and $0$ otherwise. Consider the exact
sequence~:
$$
H^1(S_x,\f O(\m d-\m l))\fl H^1(S_x,\f I_Q(\m d))\fl H^1(\til L,\f
I_{Q,\til L}(\m d))$$

The scheme $Z=P_1^{m_1-\eps_1}\cup\cdots\cup P_r^{m_r-\eps_r}$ is
clearly an $(m,a)$-configuration. Since $(\m d-\m l).\m c+1-g\geq 0$,
$h^0(C,\f O(d-1))\geq (d-1)a+1-g\geq \deg (Z\cap C)$, hence
$Z$ is a $(d-1,m,a)$-candidate (in the extended sense of remark
\ref{remark}). Moreover, $(d-1)\geq d_0(a,m)$ ; Proposition
\ref{theoalhi} shows that $\f L_x(\m d-\m l)$ is a regular system.
But $\chi(\m d-\m l)\geq\chi(\m d)-(d+1)\geq 0$, therefore 
$h^1(\m d-\m l)=0$.

Moreover, since $\til L$ is a rational curve, $|\f I_{Q,\til L}(\m
d)|$ is a regular system of degree $\m d.\m l-1$. If $\m d.\m l\geq
0$, then $h^1(\f I_{Q,\til L}(\m d))=0$, as expected. Otherwise, the
exact sequence 
$
0\fl \f O_{S_x}(\m d-\m l)\fl \f O_{S_x}(\m d)\fl \f O_{\til L}(\m
d)\fl 0$
shows that $\til L$ is a base component of $\f L_x(\m d)$. Consider
another line $L'$ of $\p^2$ containing none of the $r$ points
$P_1,\ldots,P_r$. The preceding argument is true, with $L'$ in place
of $L$, showing that the point $Q'=\til L\cap\til L'$ can not be a
base point of $\f L_x(\m d)$. This is a contradiction. 

Therefore, $h^1(\f O_{S_x}(\m d-\m l))=h^1(\f I_{Q,\til L}(\m d))=0$, and 
the first exact sequence yields $h^1(\f I_Q(\m d))=0$.

\medskip

\noindent $\bullet$ Thus $\f L_x(\m d)$ is base point free.
Bertini's theorem shows that $\D_x(\m d)$ is a smooth curve. Suppose
it is not geometrically irreducible, and denote
by $D_1,\ldots,D_l \ (l\geq 2)$ its geometrically irreducible
components  (over a bigger base
field) . Let $\m d_i\in \pic S_x$ be the class of $D_i$ $(1\leq
i\leq l)$. 

Consider two integers $1\leq i\neq j \leq l$. The curve $\D_x(\m d)$
being smooth, $D_i$ does not intersect $D_j$ and $\m d_i.\m d_j=0$.
Let $P$ be a point of $D_i$. The dimension of $\f L_x(\m d_j)$ being
positive (this system is base point free), there exists a curve
$D'_j\in \f L_x(\m d_j)$ containing $P$. But $D'_j.D_i=0$, so $D_i$ is
a component of $D'_j$, and $\f L_x(\m d_j-\m d_i)$ is effective. By
the same argument $\f L_x(\m d_i-\m d_j)$ is effective too, hence
$\m d_i=\m d_j$ for every $i\neq j$.

Thus $\m d=l\m d_1$ and $\m d^2=0$. The equality $\chi(\m d)+g(\m
d)=\m d^2+2$ yields $g(\m d)\leq 1-d$. Moreover, one can easily see
that $g(l\m d_1)=lg(\m d_1)-(l-1)$. Therefore, $lg(\m d_1)\leq l-d$
and (since $g(\m d_1)\geq 0$), $l\geq d$. The only possibility is
$l=d$. Then $D_1$ is the strict transform of a line such that
$D_1^2=0$. One may suppose that $\m d_1=(1;1)$, and $\m d=(d;d)$.  
This situation never happens since, by assumption $d\geq d_0(a,m)>m$.

\medskip

\noindent $\bullet$ We still have to prove that $\D_x(\m d)$ meets
$\til C$ transversally and that $\D_x(\m d)\cap \til C$ is
irreducible. The first point comes from Bertini's theorem, applied
to the curve $\til C$ and the restricted (base point free) linear
system $\f L_x(\m d)_{|\til C}$.

As for the second point, consider $\f G\subset S_x\times \f L(\m d)$, the
universal divisor over $\f L(\m d)$. Let $I=\f G\cap (\til C\times \f
L(\m d))$ be the intersection variety. Since $\f L(\m d)$ is base
point free, the fibers of the natural projection $I\fl \til C$ are
projective spaces of constant dimension. Therefore, $I$ is irreducible
and $\D_x(\m d)\cap \til C$, which is the generic fiber of $I\fl \f
L(\m d)$ also. \sq

\section{Proof of Theorem 1}
\label{sec.proof}

In this section we gather the preceding results to prove the 
announced theorem. Actually, the work is not made directly on the
projective plane but, rather, on the plane blown up at the $r$ points.
However, the theorem proved below clearly implies
the statement of the introduction.

\begin{theorem}
\label{theorem2}
Let $m$ be a positive integer,
$x$ the generic point of $(\p^2)^r$ and  $\m d=(d;m_1,\ldots,m_r)\in
\pic S_x$ such that $0<m_i\leq m$ ($1\leq i\leq r$). 
With the notations of \ref{theoalhi}, let 

\smallskip

$a=\bo a'(m)=\max (\bo a(m), 4 m) ; $

$\bo d'(m)= \max(d_0(a,m)+2a, a(2m+1)).$ 

\smallskip 

\noindent If $d\geq \bo d'(m)$ then
$\f L_x(\m d)$ is regular, and if $\dim \f L_x(\m d)\geq 0$, the
generic curve $\D_x(\m d)$ is geometrically irreducible, smooth, and
meets each exceptional divisor $E_i$ in $m_i$ distinct points
($1\leq i\leq r$). As a consequence, the image of $\D_x(\m d)$ on
the plane has an ordinary singularity of the prescribed multiplicity
$m_i$ at every $P_i$. 
\end{theorem}
\proof
Let $C$ be the generic curve of degree $a$ and
genus $g$.

If $\chi(\m d)\leq 0$, then the scheme $Z=P_1^{m_1}\cup\cdots\cup
P_r^{m_r}$ is a $(d,m,a)$-candidate with an empty constraint part.
By Proposition \ref{theoalhi}, $Z$ is winning, $\f L_x(\m d)$ 
is empty and the theorem is true.
If $\chi(\m d)>1$ one may consider $\chi(\m d)-1$ more
general points $P_{r+1},\ldots,P_{r+\chi(\m d)-1}$  of $\p^2$ and
study the sub-system of curves in  $\f L_x(\m d)$ passing through
these supplementary points with multiplicity at least $1$. It is 
equivalent to study the curves in $\f L(d;m_1,\ldots,m_r)$ or in 
$\f L(d;m_1,\ldots,m_r,1^{\chi(\m d)-1})$.
As a consequence, we can make the assumption that $\chi(\m d)=1$.

\medskip

\noindent $\bullet$ There exists a positive integer $s\leq r$ such that~:
\begin{eqnarray}
\label{ajuste}
-\al &= & da-m_1-\cdots -m_s +1 -g \in [-d+a-m,-d+a-1]\\
\label{bigs}
s&\geq& (2da-a^2)/(2m)
\end{eqnarray}
The second inequality follows from the first one~: 
(\ref{ajuste}) together with $m_i\leq m$ gives 
$da -ms +1-g\leq -d+a-1$ hence (since $g\leq a^2/2$),
$ms\geq (2da-a^2)/2$.
As for (\ref{ajuste}), since $0<m_i\leq m$, it is sufficient to show that 
$da-\sum_{i=1}^r m_i +1-g\leq -d+a-1$. This is a consequence of the
assumption $\chi(\m d)=1$~:
$$
\begin{array}{ll}
& da-(m_1+\cdots +m_r)+1-g\leq -d+a-1\\
\Longleftarrow_{(m_i\leq m,g>0)}& da-\frac{2}{m+1}
	\left(\frac{m_1(m_1+1)}2 +\cdots
	+\frac{m_r(m_r+1)}2\right)\leq  -d+a-1\\
\equi & da-\frac{2}{m+1}\left(\frac{(d+1)(d+2)}2-\chi(\m d)\right)\leq
	-d+a-1\\ 
\equi_{(\chi(\m d)=1)} &  d(a+1)-\frac{d(d+3)}{(m+1)}-a+1\leq 0
\end{array}
$$
which is true when $d\geq a(2m+1)$.

\noindent $\bullet$ Let $x=(P_1,\ldots,P_r)$ denote the generic point of
$(\p^2)^r$, and $y=(Q_1,\ldots,$ $Q_s,P_{s+1},\ldots,P_r)$ the generic point of $C^s\times
(\p^2)^{r-s}$. The $r$-tuple $y$ is a specialization of $x$. The class
of $\til C_y$ is $\m c=(a;1^s,0^{r-s})$. In order to
apply the Geometric Horace Lemma we are going to check the points
$\hyp 1$ to $\hyp {10}$ of lemma \ref{horgeo}.

Condition \hyp 1 is nothing but the relation (\ref{ajuste}) above. 
The assumption $d\geq a(2m+1)$ and (\ref{bigs}) yield 
$s\geq a^2\geq g$, hence \hyp 2 is true ; moreover
(\ref{bigs}) and $a\geq 4 m$ also give 
$s\geq 4 d -2a\geq d-a+m+1$, hence \hyp 3 is
true. 

As for the regularity of $\f I_{Q_1^C\cup\cdots\cup Q_{\al+1}^C}
(\m d-\m c))$, we make use of the corollary \ref{vanish}.
Since $\chi(\m d)=1$, the exact sequence \ref{exact} (section
\ref{sec.horgeo}) and \ref{ajuste} yield
$\chi(\m d-\m c)-(\al+1)=0$. The choice of $\bo a'(m)$ and
$\bo d'(m)$ gives $a\geq \max(\bo a(m),4 m)$ and
$d-a\geq \max(d_0(a,m)+a, 2am)$. Therefore, the only
remaining condition of \ref{vanish} is
$$
\begin{array}{ll}
& (\m d-\m c).\m c - \al -1 +1-g \geq 0\\
\equi & (d-a).a-\sum_{i=1}^s(m_i-1)+1-g-(\al
+1)\geq 0\\
\equi_{(\ref{ajuste})}& -2\al -a^2+s-1\geq 0 \\
\Longleftarrow_{(\ref{bigs})}& -2\al - a^2
+(2da-a^2)/(2m)-1 \geq 0\\
\Longleftarrow_{(\al \leq d-a+m)}& d\bigl(\frac a m - 2\bigr)
-a^2 \bigl( 1 +\frac 1 {2m}\bigr) +2a - 2m -1 \geq 0\\
\Longleftarrow_{(d \geq a(2m+1))}&
(2a^2 +\frac {a^2} m  -4am-2a)-a^2 - \frac{a^2}{2m}+2a-2m-1 \geq 0\\
\equi & a^2 - 4am +\frac{a^2}{2m} -2m -1 \geq 0 \\
\Longleftarrow_{(a\geq 4m)}& 8m -2m -1 \geq 0
\end{array}
$$
which is true since $m\geq 1$.

\medskip

\noindent $\bullet$
We are now left with the ``irreducibility'' and ``smoothness'' part
of lemma \ref{horgeo}. The sheaf $\f O(\m c)$ is effective if and only
if $s\leq a(a+3)/2$, which is not the case by (\ref{bigs}) ; so
$\hyp 5$ is true. The $7$-th point is empty since 
$\al\geq d-a+1>0$. 

Now, the residual system $\f L_y(\m d-\m c)$ has a ``high'' dimension~:
Precisely, the exact sequence \ref{exact} shows that
$\dim \f L_x(\m d-\m c)=\al\in [d-a+1,d-a+m]$.
Thus, the remaining condition  of
the Horace lemma can be proved with the proposition \ref{highdim}. 
By assumption $d-a\geq d_0(a,m)+1$. It is then
sufficient to prove that $(\m d-\m c).\m c+1-g\geq a$.

But the preceding computation shows that 
$(\m d-\m c).\m c+1-g\geq \al + 1\geq d-a+2 \geq _{(d\geq 3a)} a$.

\noindent $\bullet$ Let us turn now to the question of ordinary
singularities . Recalling the remark \ref{ordinary}, we only have to check
that the system $\f L_y(\m d-\m c)$ is base point free, which is 
the case by Proposition \ref{highdim}.
\sq

\end{document}